\documentclass[a4paper,11pt]{amsart}
\usepackage{amssymb}
\usepackage{amsmath}
\usepackage{mathrsfs}
\usepackage{amsmath,amssymb,amsthm,latexsym,amscd,mathrsfs}
\usepackage{indentfirst}
\usepackage{graphicx}
\usepackage{booktabs}
\usepackage{array}
\usepackage{chemarrow}
\newcommand{\ROM}[1]{\mathrm{\uppercase\expandafter{\romannumeral#1}}}
\theoremstyle{definition}

\newtheorem{thm}{Theorem}[section]

\newtheorem{rem}{Remark}[section]

\newtheorem{ack}{Acknowledgements}   
\makeatletter


\title[New examples of Willmore submanifolds in the unit sphere via isoparametric functions]{\textbf{New examples of Willmore submanifolds in the unit sphere via isoparametric functions}}
\author[Z.Z.Tang]{Zizhou Tang}\address{School of Mathematical Sciences, Laboratory of Mathematics and Complex Systems, Beijing Normal
University, Beijing 100875, China}\email{zztang@bnu.edu.cn}
\thanks {The project is partially supported by the NSFC ( No.11071018 ) and the Program for Changjiang Scholars and Innovative
Research Team in University.}
\author[W. J. Yan]{Wenjiao Yan}
\address{School of Mathematical Sciences, Laboratory of Mathematics and Complex Systems, Beijing Normal
University, Beijing 100875, China} \email{wjyan@mail.bnu.edu.cn}

\date{}
\begin{document}

\maketitle
\begin{center}
Dedicated to Professor Chiakuei Peng on his 70-th birthday.
\end{center}

\begin{abstract}
An isometric immersion $x:M^n\rightarrow S^{n+p}$ is called Willmore
if it is an extremal submanifold of the Willmore functional:
$W(x)=\int_{M^n} (S-nH^2)^{\frac{n}{2}}dv$, where $S$ is the norm
square of the second fundamental form and $H$ is the mean curvature.
Examples of Willmore submanifolds in the unit sphere are scarce in
the literature. The present paper gives a series of new examples of
Willmore submanifolds in the unit sphere via isoparametric functions
of FKM-type.
\end{abstract}

\section{Introduction}

Let $M$ be an $n$-dimensional compact submanifold immersed in an $(n + p)$-dimensional
unit sphere $S^{n+p}$. Denote by $h$ the second fundamental form of $M$, $S$
the norm square of $h$, $\vec{\textbf{H}}$ the
mean curvature vector and $H$ the mean curvature of $M$, respectively.
Based on the following range of indices:
$$1 \leq i,j,k \leq n;\quad n + 1\leq \alpha, \beta, \gamma \leq n + p;\quad 1 \leq A,B,C \leq n + p,$$
we let $\{e_A\}$ be a field of local orthonormal basis for $TS^{n+p}$
such that when restricted to $M^n$, $\{e_i\}$ is a field of local orthonormal basis for
$TM$ and $\{e_{\alpha}\}$ is a field of local orthonormal basis for the normal bundle $T^{\perp}M$.
In this way, $h$ has components $h^{\alpha}_{ij}$, and we immediately get the following expressions:
\begin{equation*}
S=\sum_{\alpha,i,j}(h^{\alpha}_{ij})^2,\quad \vec{\textbf{H}}=\sum_{\alpha}H^{\alpha}e_{\alpha},\quad H^{\alpha}=\frac{1}{n}\sum_{i}h^{\alpha}_{ii},\quad
H=|\vec{\textbf{H}}|.
\end{equation*}

$M^n$ is called a \emph{Willmore submanifold} in $S^{n+p}$ if it is
an extremal submanifold of the Willmore functional, which is
conformal invariant (\emph{cf.} \cite{Wan1}):
$$W(x)=\int_{M^n} (S-nH^2)^{\frac{n}{2}}dv.$$
Defining a non-negative function $\rho^2$ on $M$:
\begin{equation*}
\rho^2=S-nH^2,
\end{equation*}
we recall an equivalent condition for $M$ to be Willmore:

\vspace{3mm}

\noindent
\textbf{Theorem (\cite{GLW}, \cite{PW})}\,\,
{\itshape
$M^n$ is a Willmore submanifold in $S^{n+p}$ if and only if for
any $\alpha$ with $n+1\leq \alpha \leq n+p$,
\begin{eqnarray}\label{Willmore submanifold}
&&-\rho^{n-2}\sum_{i,j}\Big(R_{ij}-(n-1)\sum_{\beta}H^{\beta}h_{ij}^{\beta}\Big)\Big(h_{ij}^{\alpha}-H^{\alpha}\delta_{ij}\Big)\\
&&+(n-1)\Delta(\rho^{n-2}H^{\alpha})-\sum_{i,j}(\rho^{n-2})_{,ij}(nH^{\alpha}\delta_{ij}-h^{\alpha}_{ij})=0,\nonumber
\end{eqnarray}
where $R_{ij}$ is the Ricci tensor.}
\begin{rem}
When $n=2$, the formula (\ref{Willmore submanifold}) reduces to the following well known criterion of Willmore surfaces:
\begin{equation}\label{Willmore submanifold for n=2}
\Delta^{\perp}H^{\alpha}+\sum_{\beta,i,j}H^{\beta}h^{\alpha}_{ij}h^{\beta}_{ij}-2|\vec{\textbf{H}}|^2H^{\alpha}=0,\qquad 3\leq \alpha \leq 2+p.
\end{equation}
Clearly, every minimal surface in $S^{2+p}$ is automatically
Willmore; in other words, Willmore surfaces are a generalization of
minimal surfaces in a sphere.\hfill $\Box$
\end{rem}

It follows immediately from the theorem above that all
$n$-dimensional Einstein manifolds minimally immersed in the unit
sphere are Willmore submanifolds. However, as a matter of fact,
there do exist examples of Willmore hypersurfaces in the unit
sphere, which are neither minimal nor Einstein. For instance,
Cartan's minimal isoparametric hypersurface is Willmore but not
Einstein, while one certain hypersurface in each of Nomizu's isoparametric
families is Willmore, but neither minimal nor Einstein. In
addition, \cite{Li} characterized all isoparametric Willmore
hypersurfaces in the unit sphere.

In this paper, by taking advantage of isoparametric functions of FKM-type, we give
a series of new examples of Willmore submanifolds $M^n$ in the unit
sphere $S^{n+p}$, which are minimal but in general not Einstein (for details, see the concluding remark).

As is well known, a hypersurface $M^n$ in the unit sphere $S^{n+1}$ is \emph{isoparametric}
if it is a level hypersurface of some locally defined isoparametric function $f$ on $S^{n+1}$,
that is, a non-constant smooth function $f: S^{n+1}\rightarrow \mathbb{R}$ (\emph{cf}. \cite{Wa87}) satisfying:
\begin{equation}\label{ab}
\left\{ \begin{array}{ll}
|\nabla f|^2= b(f)\\
~~~~\triangle f~~=a(f)
\end{array}\right.
\end{equation}
where $\nabla f$ and $\triangle f$ is the gradient and Laplacian of $f$, respectively, $b$ is a smooth function on $\mathbb{R}$, and $a$ is a continuous function on $\mathbb{R}$.

\'{E}lie Cartan (\emph{cf}. \cite{Ca1,Ca2,Ca3,Ca4}) pointed out
that the level hypersurfaces $M_t:=f^{-1}(t)$ corresponding to regular
values $t$ of $f$ are parallel (which are called \emph{isoparametric hypersurfaces}) and have constant principal curvatures.
It is well known that the focal submanifolds $M_+:=f^{-1}(+1)$ and $M_-:=f^{-1}(-1)$ are minimal submanifolds
in $S^{n+1}$ (\emph{cf.} \cite{CR}). In fact, as asserted by Ge and Tang (\cite{GT}), this result still holds for focal submanifolds of isoparametric hypersurfaces in any Riemannian manifold.

We now recall the construction of isoparametric functions of
FKM-type. For a symmetric Clifford system $\{P_0,\cdots,P_m\}$ on
$\mathbb{R}^{2l}$, \emph{i.e.} $P_i$'s are symmetric matrices
satisfying $P_iP_j+P_jP_i=2\delta_{ij}I_{2l}$, Ferus, Karcher and
M\"{u}nzner (\cite{FKM}) constructed a polynomial $F$ on
$\mathbb{R}^{2l}$:
\begin{eqnarray}\label{FKM isop. poly.}
&&\qquad F:\quad \mathbb{R}^{2l}\rightarrow \mathbb{R}\nonumber\\
&&F(x) = |x|^4 - 2\displaystyle\sum_{i = 0}^{m}{\langle
P_ix,x\rangle^2}.
\end{eqnarray}
It is not difficult to verify that $f=F|_{S^{2l-1}}$ satisfies (\emph{cf.}\cite{Ce}, \cite{GTY}):
\begin{equation}\label{ab}
\left\{ \begin{array}{ll}
|\nabla f|^2= 16 (1-f^2)\\
~~~\triangle f~~=8(m_2-m_1)-4(2l+2)f.
\end{array}\right.
\end{equation}
where $m_1=m$, $m_2=l-m-1$ are the differences of the dimensions of
$M_+$ and $M_-$ compared with that of the isoparametric
hypersurface, respectively, $\nabla f$ and $\triangle f$ are the
gradient and Laplacian of $f$ on $S^{2l-1}$, respectively. Thus by
definition, $f$ is an \emph{isoparametric function} on $S^{2l-1}$,
which is called the \emph{FKM-type isoparametric polynomial} on
$S^{2l-1}$. Correspondingly, the focal submanifolds
$M_+:=f^{-1}(+1)$ can be expressed as an algebraic set:
$$M_+=\{x\in S^{2l-1}~|~\langle P_0x,x\rangle=...=\langle P_mx,x\rangle=0~\},$$
which will be exactly our new examples of Willmore submanifolds in
$S^{2l-1}$.

The main result of the present paper is stated as follows:

\begin{thm}\label{prop1}
{\itshape
The focal submanifolds $M_+$ of the FKM-type isoparametric polynomials on $S^{2l-1}$ are Willmore submanifolds in $S^{2l-1}$.}
\end{thm}

\section{Proof of Theorem \ref{prop1}}

We start with a description of the normal space of the focal
submanifold $M_+$ following \cite{FKM}. Since $M_+=\{~x\in
S^{2l-1}~|~\langle P_0x,x\rangle=...=\langle P_mx,x\rangle=0~\}$, it
follows without much difficulty that $dim~M_+=2l-1-(m+1)=2l-m-2$,
and as pointed out by \cite{FKM}, the normal space at $x\in M_+$ is
\begin{equation}
T^{\perp}_xM_+=\{~\mathcal{P}x~|~\mathcal{P}\in\mathbb{R}\Sigma (P_0,...,P_m)~\},
\end{equation}
where $\Sigma (P_0,...,P_m)$ is the unit sphere in $Span\{P_0,...,P_m\}$, which is called \emph{the Clifford sphere}
determined by the system $\{P_0,...,P_m\}$.

For a normal vector $\xi_{\alpha}=P_{\alpha}x$, $\alpha=0,...,m$, let $\langle A_{\alpha}X, Y\rangle=\langle h(X, Y), \xi_{\alpha}\rangle$,
$\forall X,Y \in T_xM_+$, where $A_{\alpha}$ is the shape operator corresponding to $\xi_{\alpha}$; in other words, we have $h(X, Y)=\sum_{\alpha}\langle A_{\alpha}X, Y\rangle\xi_{\alpha}$.

Next, from Gauss equation, we derive that
\begin{equation}\label{Gauss}
K(X, Y)= 1+\sum_{\alpha=0}^m\Big\{ \langle A_{\alpha}X,
X\rangle\langle A_{\alpha}Y, Y\rangle - \langle A_{\alpha}X,
Y\rangle^2\Big\},
\end{equation}
where $K$ is the sectional curvature in $M_+$. Furthermore, in
virtue of the fact $A_{\alpha}X=-(P_{\alpha}X)^T$ (\emph{cf}.
\cite{FKM}), the tangential component of $-P_{\alpha}X$, the formula
(\ref{Gauss}) becomes:
\begin{equation}\label{sectional curvature}
K(X, Y)= 1+\sum_{\alpha=0}^m\Big\{\langle P_{\alpha}X, X\rangle\langle P_{\alpha}Y, Y\rangle - \langle P_{\alpha}X, Y\rangle^2\Big\}
\end{equation}

Let $\{X=e_1$, $e_2,...,e_{2l-m-2}\}$ be an orthonormal basis for $T_xM_+$, and $\{P_0x,...,P_mx\}$ an orthonormal basis
for $T^{\perp}_xM_+$. Then the Ricci curvature for $X$ is
\begin{eqnarray}\label{Ricci}
Ricci (X)&=&\sum_{i=2}^{2l-m-2}K(X, e_i)\\
&=& 2l-m-3 + \sum_{i=2}^{2l-m-2}\sum_{\alpha=0}^{m}\Big\{\langle P_{\alpha}X, X\rangle\langle P_{\alpha}e_i, e_i\rangle - \langle P_{\alpha}X, e_i\rangle^2\Big\}\nonumber\\
&=& 2l-m-3-\sum_{\alpha=0}^{m}\Big\{\langle P_{\alpha}X, X\rangle^2 + \sum_{i=2}^{2l-m-2}\langle P_{\alpha}X, e_i\rangle^2 \Big\}\nonumber\\
&=& 2(l-m-2) + 2\sum_{\alpha, \beta =0, \alpha < \beta}^{m}\langle X, P_{\alpha}P_{\beta}x\rangle^2,\nonumber
\end{eqnarray}
where the third equality is derived from the fact that for any
$\alpha \in \{0,...,m\}$, $P_{\alpha}$ is trace free, and
\begin{eqnarray}\label{vanishing trace}
trace~P_{\alpha}&=&\langle P_{\alpha}X, X\rangle + \sum_{i=2}^{2l-m-2}\langle P_{\alpha}e_i, e_i\rangle + \sum_{\beta=0}^{m}\langle P_{\alpha}P_{\beta}x, P_{\beta}x\rangle + \langle P_{\alpha}x, x\rangle\\
&=& \langle P_{\alpha}X, X\rangle + \sum_{i=2}^{2l-m-2}\langle P_{\alpha}e_i, e_i\rangle,\nonumber
\end{eqnarray}
and the fourth equality is obtained by the relation that for any
$\alpha \in \{0,...,m\}$,
\begin{eqnarray}\label{length}
1 &=& |P_{\alpha}X|^2\\
&=& \langle P_{\alpha}X, X\rangle^2 + \sum_{i=2}^{2l-m-2}\langle P_{\alpha}X, e_i\rangle^2 + \sum_{\beta=0}^{m}\langle P_{\alpha}X, P_{\beta}x\rangle^2 + \langle P_{\alpha}X, x\rangle^2 \nonumber\\
&=& \langle P_{\alpha}X, X\rangle^2 + \sum_{i=2}^{2l-m-2}\langle P_{\alpha}X, e_i\rangle^2 + \sum_{\beta=0}^{m}\langle X, P_{\alpha}P_{\beta}x\rangle^2.\nonumber
\end{eqnarray}
Moreover, the the norm square of the second fundamental form of
$M_+$ can be calculated immediately

\begin{eqnarray}\label{S}
S &=& \sum_{i,j=1}^{2l-m-2}\sum_{\alpha=0}^{m} \langle A_{\alpha}e_i, e_j\rangle^2\\
&=& \sum_{i,j=1}^{2l-m-2}\sum_{\alpha=0}^{m} \langle P_{\alpha}e_i, e_j\rangle^2\nonumber\\
&=& \sum_{i=1}^{2l-m-2}\sum_{\alpha=0}^{m}\Big\{ 1- \sum_{\beta=0}^{m}\langle P_{\alpha}e_i, P_{\beta}x\rangle^2 \Big\}\nonumber\\
&=& (2l-m-2)(m+1)- 2\sum_{\alpha, \beta =0, \alpha < \beta}^{m}|P_{\alpha}P_{\beta}x|^2\nonumber\\
&=& 2(l-m-1)(m+1).\nonumber
\end{eqnarray}
\begin{rem}
In fact, for every isoparametric hypersurface with four distinct
principal curvatures in the unit sphere, the norm square $S$ of the
second fundamental forms of both focal submanifolds are constant.
The proof is a consequence of a result of M\"{u}nzner (cf. \cite{CR}
page 248).
\end{rem}

As a direct result of (\ref{S}), $\rho^2=S-nH^2=S$ is a constant on
$M_+$. Together with the minimality of $M_+$, \emph{i.e.}
$H^{\alpha}=0$ for any $0 \leq \alpha \leq m$, the criterion
(\ref{Willmore submanifold}) for Willmore reduces to
\begin{equation}\label{reduced Willmore}
for ~any ~ 0\leq \alpha \leq m, ~~\sum_{i,j=1}^{2l-m-2} R_{ij}h^{\alpha}_{ij}=0.
\end{equation}

In order to prove that (\ref{reduced Willmore}) holds for $M_+$, we
first note that for any $\mathcal{P}\in \Sigma(P_0,...,P_m)$, the
eigenvalues of $\mathcal{P}$ must be $\pm 1$, since $\mathcal{P}^2=
I$. Observe that since $trace ~\mathcal{P}=0$, $E_+$ and $E_-$, the
eigenspaces of $\mathcal{P}$ for the eigenvalues $+1$ and $-1$, have
the same dimension $l$, so that $\mathbb{R}^{2l}=E_+\oplus E_-$.
Then we can compute the principal curvatures of the shape operator
$A_{\xi}$, with respect to a unit normal $\xi=\mathcal{P}x$ as
follows.

\vspace{2mm}

\noindent \textbf{Lemma (\cite{Ce}, \cite{FKM})}\,\, {\itshape Let
$x$ be a point in the focal submanifold $M_+$, and let
$\xi=\mathcal{P}x$ be a unit normal vector to $M_+$ at $x$, where
$\mathcal{P}\in\Sigma(P_0,...,P_m):=\Sigma$. Then the shape operator
$A_{\xi}$ has principal curvatures $0$, $1$, $-1$ with corresponding
principal spaces $T_0(\xi)$, $T_1(\xi)$ and $T_{-1}(\xi)$ as
follows:
\begin{eqnarray}\label{eigenspaces}
T_0(\xi)&=& \{ ~QPx~|~Q\in \Sigma ,~ \langle Q, P\rangle =0 ~\},\\
T_1(\xi)&=& \{ ~X\in E_-~|~X\cdot Qx =0,~ \forall Q\in \Sigma ~\} = E_-\cap T_xM_+,\nonumber\\
T_{-1}(\xi)&=& \{ ~X\in E_+~|~X\cdot Qx =0,~ \forall Q\in \Sigma ~\} = E_+\cap T_xM_+.\nonumber
\end{eqnarray}
Furthermore,
$$dim T_0(\xi)=m,\quad dim T_1(\xi)=dim T_{-1}(\xi)=l-m-1.$$
}

In this way, $T_xM_+$ can be decomposed as a direct sum
$T_xM_+=T_0(\xi)\oplus T_1(\xi) \oplus T_{-1}(\xi)$. Let's choose an
orthonormal basis $\{u_1,...,u_m,v_1,...,v_{l-m-1},
w_1,...,w_{l-m-1}\}$ on $T_xM_+$ with $u_1,...,u_m\in T_0(\xi)$,
$v_1,...,v_{l-m-1}\in T_1(\xi)$ and $w_1,...,w_{l-m-1}\in
T_{-1}(\xi)$. Then it is easy to see
$$\sum_{i,j=1}^{2l-m-2} R_{ij}h^{\alpha}_{ij}=\sum_{i=1}^{l-m-1}\{Ricci (v_i)-Ricci (w_i)\}.$$

Thus it follows from (\ref{reduced Willmore}) that
\begin{equation}\label{equivalence for v w}
M_+ ~is ~Willmore~\Leftrightarrow~\sum_{i=1}^{l-m-1}Ricci (v_i)=
\sum_{i=1}^{l-m-1}Ricci (w_i).
\end{equation}

However, by (\ref{Ricci}), we have
\begin{eqnarray}\label{T1 part}
&&\sum_{i=1}^{l-m-1}Ricci (v_i)=\sum_{i=1}^{l-m-1}Ricci (w_i)\nonumber\\
&\Leftrightarrow& \sum_{i=1}^{l-m-1}\sum_{\alpha, \beta =0}^m\langle P_{\alpha}v_i, P_{\beta}x\rangle^2 = \sum_{i=1}^{l-m-1}\sum_{\alpha, \beta =0}^m\langle P_{\alpha}w_i, P_{\beta}x\rangle^2\\
&\Leftrightarrow& \sum_{\alpha, \beta=0, \alpha \neq
\beta}^m|(P_{\alpha}P_{\beta}x)^{T_1}|^2 = \sum_{\alpha, \beta=0,
\alpha \neq \beta}^m|(P_{\alpha}P_{\beta}x)^{T_{-1}}|^2\nonumber
\end{eqnarray}

Finally, we notice that when $m=1$, the rightmost side of (\ref{T1
part}) automatically holds. In order to complete the proof of the
theorem, we need the following observation.

Given any $\mathcal{P}\in \Sigma(P_0,...,P_m)$, we can assume
$\mathcal{P}=P_0$ without loss of generality ( choose a suitable
orthogonal transformation on $\mathbb{R}\Sigma(P_0,...,P_m)$
preserving the geometric equivalence class ). Under this assumption,
we get a consequence that for $X\in T_1(\xi)$, $P_0X=-X$, and for
$Y\in T_{-1}(\xi)$, $P_0Y=Y$. In fact, by $(P_0X)^T=-A_{\xi}X=-X$,
we can decompose $P_0X$ as $P_0X=-X+(P_0X)^N$, where $(P_0X)^N$ is
the normal part of $P_0X$. Then from the identities $|P_0X|=|X|=1$,
we derive that $P_0X=-X$. The other proof with respect to $Y$ is
analogous.

Now with this preparation, it suffices to prove the rightmost side of (\ref{T1 part}) in the following
two cases.

Case A: $m=2$. In this case, $P_0P_1x$, $P_0P_2x$ $\in T_0(\xi)$, we
just need to deal with $P_1P_2x$. Since $\langle P_1P_2x,
P_{\alpha}x\rangle =0$ for any $\alpha=0,1,2,$ and $\langle P_1P_2x,
x\rangle =0$, we know $P_1P_2x\in T_xM_+$.

Decompose $P_1P_2x$ as $P_1P_2x=U+V+W\in T_0(\xi)\oplus
T_1(\xi)\oplus T_{-1}(\xi)$. In virtue of the observation, we have
$P_0P_1P_2x=P_0U-V+W=-V+W$. In fact, according to
(\ref{eigenspaces}), we can write $U=Q_0P_0x$ with $Q_0\in \Sigma$,
so that $P_0U=-Q_0x\in T^{\perp}_xM_+$. On the other hand, since
$\langle P_0P_1P_2x, P_{\alpha}x\rangle =0$ for any $\alpha=0,1,2$,
we know that $P_0P_1P_2x\in T_xM_+$. Thus $P_0U=0$.

From $\langle P_1P_2x, P_0P_1P_2x\rangle =0$, it follows directly that $\langle -V+W, U+V+W\rangle=-|V|^2+|W|^2=0$,
or equivalently, $|(P_1P_2x)^{T_1}|=|(P_1P_2x)^{T_{-1}}|$.
\vspace{2mm}

Case B: $m>2$. In this case, $P_0P_{\alpha}x\in T_0(\xi)$ for any
$\alpha = 1,...,m$, so we just need to deal with
$P_{\alpha}P_{\beta}x$ for $\alpha, \beta >0$. Observing that for
$\alpha \neq \beta$, $\langle P_{\alpha}P_{\beta}x,
P_{\gamma}x\rangle =0$ for any $\gamma=0,...,m$, and $\langle
P_{\alpha}P_{\beta}x, x\rangle =0$, we have $P_{\alpha}P_{\beta}x\in
T_xM_+$.

Decomposing $P_{\alpha}P_{\beta}x$ as $P_{\alpha}P_{\beta}x=U+V+W\in
T_0(\xi)\oplus T_1(\xi)\oplus T_{-1}(\xi)$, we get
$P_0P_{\alpha}P_{\beta}x=P_0U-V+W$, and
\begin{equation}\label{plus,minus}
\left\{ \begin{array}{ll}
P_{\alpha}P_{\beta}x + P_0P_{\alpha}P_{\beta}x = U+P_0U+2W\\
P_{\alpha}P_{\beta}x - P_0P_{\alpha}P_{\beta}x = U-P_0U+2V.
\end{array}\right.
\end{equation}
Taking norm square on both sides of the two equalities above, since
$\langle P_0P_{\alpha}P_{\beta}x, P_{\alpha}P_{\beta}x\rangle=0$ for
any $\alpha, \beta >0$, we arrive at
\begin{equation}\label{|V|=|W|}
\left\{ \begin{array}{ll}
2 = |U|^2+|P_0U|^2+4|W|^2\\
2 = |U|^2+|P_0U|^2+4|V|^2,
\end{array}\right.
\end{equation}
which implies $|V|^2=|W|^2$, \emph{i.e.}
$|(P_{\alpha}P_{\beta}x)^{T_1}|=|(P_{\alpha}P_{\beta}x)^{T_{-1}}|$.

The proof of our main theorem is now complete!
\vspace{5mm}

\noindent \textbf{Concluding Remark.}\,\, { In general, the focal
submanifold $M_+$ of the FKM-type isoparametric function is not
Einstein. As is well known, for each positive integer $k$, there
exists one Clifford system $\{P_0,...,P_m\}$ on $\mathbb{R}^{2l}$
with $l=k\delta(m)$, where $\delta(m)$ has the following values:

\begin{center}
\begin{tabular}{|c|c|c|c|c|c|c|c|c|c|}
\hline
$m$ & 1 & 2 & 3 & 4 & 5 & 6 & 7 & 8 & $\cdots$ $m$+8 \\
\hline
$\delta(m)$ & 1 & 2 & 4 & 4 & 8 & 8 & 8 & 8 & ~16$\delta(m)$\\
\hline
\end{tabular}
\end{center}

Since each focal manifold of the FKM-type isoparametric function has
reduced dimension compared with the isoparametric hypersurface, we
have $m_1=m\geq 1$ and $m_2=l-m-1\geq 1.$ On these two conditions,
with a fundamental but complicated argument (the details are
omitted), we observe that if $m=1, 2, 3$ or $m>10$, then
$4l>m^2+3m+4$, which implies that
$$dim~M_+>\frac{1}{2}m(m+1)\geq dim~Span\{P_{\alpha}P_{\beta}x~|~ \alpha,\beta=0,1,...,m, \alpha<\beta\}.$$
At last, combining with (\ref{Ricci}), we can conclude that the
Ricci curvature of $M_+$ is not constant in these cases, namely,
$M_+$ is not an Einstein manifold ! }

\begin{ack}
The authors would like to thank Professor Haizhong Li for his providing the pdf file of the reference
\cite{PW} of Pedit and Willmore.
\end{ack}

\end{document}